\documentclass{amsart}
\usepackage{graphicx}
\usepackage{tikz-cd}
\usepackage{xcolor}
\usepackage{comment}
\usepackage[T1]{fontenc}
\usepackage{xypic}

\usepackage{amsmath,amsthm,amssymb,amscd,enumerate}
\usepackage{enumitem}

\usetikzlibrary{decorations.pathmorphing}

\newtheorem{theorem}{Theorem}[section]
\newtheorem{proposition}[theorem]{Proposition}

\newtheorem{corollary}[theorem]{Corollary}

\theoremstyle{definition}

\theoremstyle{remark}

\newtheorem{remark}[theorem]{Remark}

\numberwithin{equation}{section}

\allowdisplaybreaks[1]

\newcommand{\QQ}{\mathbb{Q}}
\newcommand{\RR}{\mathbb{R}}
\newcommand{\ZZ}{\mathbb{Z}}

\newcommand{\QQQ}{\mathrm{Q}}

\newcommand{\HHH}{\mathrm{H}}

\DeclareMathOperator{\scl}{\mathrm{scl}}

\newcommand{\Gam}{\Gamma}

\newcommand{\WW}{\mathrm{W}}
\newcommand{\VV}{\mathrm{V}}

\newcommand{\DD}{D}

\renewcommand{\Im}{\mathrm{Im}}

\newcommand{\Homeo}{\mathrm{Homeo}}

\newcommand{\Sp}{\mathrm{Sp}}
\newcommand{\Mod}{\mathrm{Mod}}

\makeatletter
\@namedef{subjclassname@2020}{%
  \textup{2020} Mathematics Subject Classification}
\makeatother


\title[Extendability problem and mapping class groups]{Note on the extendability problem for invariant quasimorphisms for mapping class groups}

\author[M. Kawasaki]{Morimichi Kawasaki}
\address[Morimichi Kawasaki]{Department of Mathematics, Faculty of Science, Hokkaido University, North 10, West 8, Kita-ku, Sapporo, Hokkaido 060-0810, Japan}
\email{kawasaki@math.sci.hokudai.ac.jp}

\author[M. Kimura]{Mitsuaki Kimura}
\address[Mitsuaki Kimura]{Department of Mathematics, Osaka Dental University,
 8-1 Kuzuha-hanazono-cho, Hirakata, Osaka 573-1121, Japan}
\email{kimura-m@cc.osaka-dent.ac.jp}

\author[S. Maruyama]{Shuhei Maruyama}
\address[Shuhei Maruyama]{School of Mathematics and Physics, College of Science and Engineering, Kanazawa University, Kakuma-machi, Kanazawa, Ishikawa, 920-1192, Japan}
\email{smaruyama@se.kanazawa-u.ac.jp}

\author[T. Matsushita]{Takahiro Matsushita}
\address[Takahiro Matsushita]{Department of Mathematical Sciences, Shinshu University, Matsumoto, Nagano 390-8621, Japan}
\email{matsushita@shinshu-u.ac.jp}

\author[M. Mimura]{Masato Mimura}
\address[Masato Mimura]{Mathematical Institute, Tohoku University, 6-3, Aramaki Aza-Aoba, Aoba-ku, Sendai 980-8578, Japan}
\email{m.masato.mimura.m@tohoku.ac.jp}

\begin{document}

\begin{abstract}
In this note, we study the extendability problem of invariant quasimorphisms for short exact sequences related to mapping class groups of surfaces, possibly with one boundary component or one marked point. 
Using a recent result of Blatz, we prove that if the genus is at least two, then every homogeneous quasimorphism on the Torelli group that is invariant by the conjugation action by the mapping class group is uniquely extendable to the whole mapping class group. 
Contrastingly, for the surface with one marked point whose genus is at least three, we show the existence of a homogeneous Torelli-group-invariant quasimorphism on the Johnson kernel that is not extendable. 
\end{abstract}

\maketitle

\section{Main results}\label{section=intro}
Let $G$ be a group and $N$ a normal subgroup of $G$, and set $\Gam=G/N$. In other words, we consider a short exact sequence of groups
\begin{equation}\label{eq=exactseq}
1\longrightarrow N\stackrel{i}{\longrightarrow} G\stackrel{p}{\longrightarrow} \Gam\longrightarrow 1,
\end{equation}
where $i\colon N\hookrightarrow G$ and $p\colon G\twoheadrightarrow \Gam$ denote the inclusion map and the group quotient map, respectively. 
A \emph{quasimorphism} on $G$ is a function $\phi\colon G\to \RR$ with
\[
\DD(\phi):=\sup\{|\phi(g_1g_2)-\phi(g_1)-\phi(g_2)|\;|\;g_1,g_2\in G\}<\infty;
\]
such $\phi$ is said to be \emph{homogeneous} if $\phi(g^k)=k\phi(g)$ for every $g\in G$ and $k\in \ZZ$. We write $\QQQ(G)$ for the $\RR$-linear space of homogeneous quasimorphisms on $G$. We say a function $\mu\colon N\to \RR$ is \emph{$G$-invariant} if $\mu(gxg^{-1})=\mu(x)$ for every $x\in N$ and $g\in G$. We write $\QQQ(N)^G$ for the $\RR$-linear space of homogeneous $G$-invariant quasimorphisms on $N$. A key observation here is $\QQQ(G)=\QQQ(G)^G$, that is, every homogeneous quasimorphism on $G$ is automatically $G$-invariant; see \cite[Lemma~1.4]{KKMMMsurvey}. Thus, $i\colon N\hookrightarrow G$ induces an $\RR$-linear map
\[
i^{\ast}\colon \QQQ(G)=\QQQ(G)^G\to \QQQ(N)^G;\quad \phi\mapsto i^{\ast}\phi=\phi\circ i(=\phi|_N).
\]
Elements in $i^{\ast}\QQQ(G)$ are called \emph{extendable} (invariant) quasimorphisms. For such quasimorphisms, there is no mystery for the $G$-invariance. In \cite{KKMMM}, the present authors have introduced the $\RR$-linear space
\[
\VV(G,N)=\QQQ(N)^G/i^{\ast}\QQQ(G),
\]
which may be regarded as the space of \emph{non-extendable} (homogeneous) invariant quasimorphisms; the symbol $\VV(G,N)$ has been given in  \cite{KKMMMsurvey}. The \emph{extendability problem} asks for a given $\mu\in \QQQ(N)^G$, whether $\mu$ represents the zero-element in $\VV(G,N)$. 
In other words, it asks whether $\mu$ is extendable to $G$. We refer the reader to the survey \cite{KKMMMsurvey} for more backgrounds and examples of invariant quasimorphisms, including this extendability problem.

Now, let $g$ be an integer at least two. Let $\Sigma_g$ be  a connected oriented closed surface of genus $g$. Let  $\Sigma_g^1$ and $\Sigma_{g,1}$ be a connected oriented surface of genus $g$ with one boundary component, and that with one marked point, respectively. For  $\Sigma$ being either of $\Sigma_g$, $\Sigma_g^1$ and $\Sigma_{g,1}$, we have the following exact sequence
\begin{equation}\label{eq=exactMod}
1\longrightarrow \mathcal{I}(\Sigma) \longrightarrow \Mod(\Sigma)\longrightarrow \Sp(2g,\ZZ)\longrightarrow 1,
\end{equation}
where $\Mod(\Sigma)$ and $\mathcal{I}(\Sigma)$ denote the mapping class group and the Torelli group of $\Sigma$, respectively. For $\Sigma=\Sigma_{g,1}$ for $g\geq 3$, we  pass one step down the Johnson filtration, and consider the following exact sequence
\begin{equation}\label{eq=exactJohnson}
1\longrightarrow \mathcal{K}_{g,1} \longrightarrow \mathcal{I}_{g,1}\stackrel{\tau}{\longrightarrow} \Lambda^3H \longrightarrow 1
\end{equation}
as well. Here, $\mathcal{I}_{g,1}$ denotes the Torelli group $\mathcal{I}(\Sigma_{g,1})$ of $\Sigma_{g,1}$, $\mathcal{K}_{g,1}$ denotes the Johnson kernel of $\Sigma_{g,1}$, and $H=\HHH_1(\Sigma_{g,1};\ZZ)$. This map $\tau$ is called the Johnson homomorphism. 

In this note, we study the extendability problems associated with the short exact sequences \eqref{eq=exactMod} and \eqref{eq=exactJohnson}. One motivation comes from the work of Cochran, Harvey and Horn \cite{CHH}. Although the Meyer signature cocycle vanishes on the Torelli group, they refined the signature construction using twisted signature defects associated with unitary representations of $\HHH_1(\Sigma;\ZZ)$. The resulting higher-order $\rho$-invariants yield, after homogenization, a family of elements of $\QQQ(\mathcal{I}(\Sigma))^{\Mod(\Sigma)}$. It is natural to ask whether the members of this family extend to $\Mod(\Sigma)$. We answer this question \emph{in the affirmative}. Our proof uses no special properties of the Cochran--Harvey--Horn quasimorphisms; it relies only on properties of the ordinary and bounded cohomology of the quotient $\Sp(2g,\ZZ)$ in \eqref{eq=exactMod}. Here is our first main result.

\begin{theorem}[vanishing of $\VV(G,N)$ for $\mathrm{Sp}(2g,\ZZ)$-quotients]\label{thm=main}
Let $\Gam$ be either of the following two groups.
\begin{enumerate}[label=\textup{(\alph*)}]
  \item $\mathrm{Sp}(2g,\ZZ)$ for $g\geq 2$.
  \item a finite index subgroup of $\mathrm{Sp}(2g,\ZZ)$ for $g\geq 3$.
\end{enumerate}
Let $(G,N)$ be a group pair that fits in the short exact sequence \eqref{eq=exactseq}. Then, the map $i^{\ast}\colon \QQQ(G)\to \QQQ(N)^G$ is an isomorphism. In other words, every homogeneous $G$-invariant quasimorphism $\mu$ on $N$ is uniquely extendable to $G$, meaning that, there uniquely exists $\phi\in\QQQ(G)$ such that $\phi|_N=\mu$. In particular, we have $\VV(G,N)=0$.
\end{theorem}

We note that option (b) covers option (a) except the group $\mathrm{Sp}(4,\ZZ)$. This subtlety for the case of $g=2$ comes from the Borel stable range for second cohomology of $\mathrm{Sp}(2g,\mathbb{Z})$; compare with the proof of Theorem~\ref{thm=Spc^2} below.
Theorem~\ref{thm=main} for (a) immediately provides the desired answer as follows.

\begin{corollary}[extendability of the Cochran--Harvey--Horn quasimorphisms]\label{cor=CHH}
Let $g$ be an integer at least two, and $\Sigma=\Sigma_g^1$. Then, every invariant quasimorphism $\mu\in \QQQ(\mathcal{I}(\Sigma))^{\Mod(\Sigma)}$ constructed by Cochran, Harvey and Horn is uniquely extendable to $\Mod(\Sigma)$.
\end{corollary}

More generally, Theorem~\ref{thm=main} provides the affirmative answer to the extendability problem in the setting of \eqref{eq=exactMod}, or its finite index analog, as follows.

\begin{corollary}\label{cor=Problem8.18}
Let $g$ be an integer at least two, and $\Sigma$ be one of $\Sigma_g$, $\Sigma_g^1$ and $\Sigma_{g,1}$. 
\begin{enumerate}[label=\textup{(\arabic*)}]
  \item Let $G=\Mod(\Sigma)$ and $N=\mathcal{I}(\Sigma)$. Then, every $\mu\in \QQQ(N)^G$ is uniquely extendable to $G$.
  \item Assume that $g\geq 3$. Let $G$ be a finite index subgroup of $\Mod(\Sigma)$ and $N=\mathcal{I}(\Sigma)\cap G$. Then, every $\mu\in \QQQ(N)^G$ is uniquely extendable to $G$.
\end{enumerate}
\end{corollary}

Corollary~\ref{cor=Problem8.18}~(2) for $\Sigma=\Sigma_{g}$ answers \cite[Problem~8.18]{KKMMM} in the affirmative. 

The key to our proof of Theorem~\ref{thm=main} is the following computation of the continuous bounded cohomology of $\Sp(2g,\RR)$ by Blatz in his Ph.D. Thesis \cite{Blatz}. Here, we note the following well-known fact: $\HHH^3_{c}(\Sp(2g,\RR);\RR)=0$ for $g\geq 2$; see \cite{BorelWallach} for instance. In addition, despite the fact that the main theorem of \cite{Blatz} is stated in terms of $\HHH^{\bullet}_{mb}(G)$, it equals $\HHH^{\bullet}_{cb}(G;\RR)$; see \cite[p.8]{Blatz}.

\begin{theorem}[{Blatz \cite[Main Theorem]{Blatz}}]\label{thm=Blatz}
For every integer $g$ at least two, the third comparison map from continuous bounded cohomology to continuous cohomology 
\[
c^3_{\Sp(2g,\RR)}\colon \HHH^3_{cb}(\Sp(2g,\RR);\RR)\to \HHH^3_{c}(\Sp(2g,\RR);\RR)
\]
is an isomorphism. In particular, we have $\HHH^3_{cb}(\Sp(2g,\RR);\RR)=0$.
\end{theorem}

More precisely, Theorem~\ref{thm=Blatz} implies that third bounded cohomology (with trivial real coefficients) $\HHH^3_b(\Gam;\RR)$ for a lattice $\Gam$ in $\Sp(2g,\RR)$ vanishes (see Theorem~\ref{thm=SpH^3_b} below), which suffices to prove Theorem~\ref{thm=main}. In \cite[Theorem~8.14]{KKMMM} the authors obtained some extendability result for quasi-equivariant vector-valued quasimorphisms. There, they employed the vanishing of third bounded cohomology $\HHH^3_b(\Gam;\mathcal{H})$ with unitary coefficients, where the invariant part of $\mathcal{H}$ is the zero space. The vanishing of $\HHH^3_b(\Gam;\RR)$ was the missing part towards Theorem~\ref{thm=main} before the aforementioned result of Blatz, and that is why \cite[Problem~8.18]{KKMMM} was raised as an open problem at that moment.

In contrast to Theorem~\ref{thm=main}, we have the \emph{negative} answer to the extendability problem in the setting of \eqref{eq=exactJohnson}. In fact, we have the non-vanishing result (Theorem~\ref{thm=main_Johnson} below) of $\WW(G,N)$, which is a quotient $\RR$-linear space of $\VV(G,N)$ defined as follows. For the short exact sequence \eqref{eq=exactseq}, the $\RR$-linear space $\WW(G,N)$ is defined by
\[
\WW(G,N)=\QQQ(N)^G/\big(\HHH^1(N;\RR)^G+i^{\ast}\QQQ(G)\big).
\]
The space $\WW(G,N)$ measures extendability with ignoring the ambiguity coming from invariant homomorphisms. This space serves as \emph{the} obstruction to the bi-Lipschitz equivalence between the stable commutator length $\scl_G$ and the stable mixed commutator length $\scl_{G,N}$ on the mixed commutator subgroup $[G,N]$; we refer the reader to \cite{KKMMM_scl1} for more details. Our second main result is the following.

\begin{theorem}\label{thm=main_Johnson}
Let $g$ be an integer at least three. Then, we have $\WW(\mathcal{I}_{g,1},\mathcal{K}_{g,1})\ne 0$. In particular, $\VV(\mathcal{I}_{g,1},\mathcal{K}_{g,1})\ne 0$.
\end{theorem}

We prove Theorem~\ref{thm=main} and Theorem~\ref{thm=main_Johnson} in Section~\ref{section=proof1} and Section~\ref{section=proof2}, respectively.

\section{Proof of Theorem~\ref{thm=main}}\label{section=proof1}
The following criterion for the vanishing of $\VV(G,N)$ is the key to the proof of Theorem~\ref{thm=main}. We refer the reader to \cite{Monodbook,Frigerio} for treatises on bounded cohomology.

\begin{proposition}\label{prop=vanishingV}
Let $\Gam$ be a group. Assume that the following hold.
\begin{enumerate}[label=\textup{(\arabic*)}]
  \item The second comparison map $c^2_{\Gam}\colon \HHH^2_b(\Gam;\RR)\to \HHH^2(\Gam;\RR)$ is surjective.
  \item The third bounded cohomology $\HHH^3_b(\Gam;\RR)$ vanishes.
\end{enumerate}
Then, for every group pair $(G,N)$ that fits in the short exact sequence \eqref{eq=exactseq}, we have $\VV(G,N)=0$. Moreover, if $\QQQ(\Gam)=0$, then the map $i^{\ast}\colon \QQQ(G)\to \QQQ(N)^G$ is an isomorphism.
\end{proposition}

For the proof of Proposition~\ref{prop=vanishingV}, we use \emph{relative bounded cohomology} $\HHH^{\bullet}_{/b}(K)$ for a group $K$, which is homology of the cochain complex $(C^{\bullet}_{/b}(K;\RR))=(C^{\bullet}(K;\RR)/C^{\bullet}_b(K;\RR))$; see \cite{KKMMM} for details. The short exact sequence of cochain complexes
\[
0\longrightarrow C^{\bullet}_b(K;\RR)\longrightarrow C^{\bullet}(K;\RR)\longrightarrow C^{\bullet}_{/b}(K;\RR)\longrightarrow 0
\]
provides the following exact sequence
\begin{equation}\label{eq=/bcohomology}
\HHH^1_{/b}(K;\RR)\longrightarrow \HHH^2_b(K;\RR)\stackrel{c^2_K}{\longrightarrow}\HHH^2(K;\RR)\longrightarrow \HHH^2_{/b}(K;\RR)\longrightarrow \HHH^3_{b}(K;\RR).
\end{equation}

\begin{proof}[Proof of Proposition~\textup{\ref{prop=vanishingV}}]
For \eqref{eq=exactseq}, Theorem 1.5 of \cite{KKMMM} states that there exists an exact sequence
\[
0\longrightarrow \HHH^1_{/b}(\Gam;\RR)\stackrel{p^{\ast}}{\longrightarrow} \HHH^1_{/b}(G;\RR) \stackrel{i^{\ast}}{\longrightarrow} \HHH^1_{/b}(N;\RR)^G \longrightarrow \HHH^2_{/b}(\Gam;\RR).
\]
Also, it is known that $\HHH^1_{/b}(K;\RR)\cong \QQQ(K)$ for a group $K$; compare with \cite[Theorem~2.50]{Calegari}. Thus, we have the following exact sequence 
\begin{equation}\label{eq=mainexact}
0\longrightarrow \QQQ(\Gam)\stackrel{p^{\ast}}{\longrightarrow} \QQQ(G) \stackrel{i^{\ast}}{\longrightarrow} \QQQ(N)^G \longrightarrow \HHH^2_{/b}(\Gam;\RR);
\end{equation}
we refer the reader to \cite[Theorem~1.5 and Remark~1.6]{KKMMM} for more details. By assumptions (1) and (2) of Proposition~\ref{prop=vanishingV}, the exact sequence \eqref{eq=/bcohomology} for $K=\Gamma$ implies that $\HHH^2_{/b}(\Gam;\RR)=0$. Therefore, by \eqref{eq=mainexact} the map $i^{\ast}\colon \QQQ(G)\to \QQQ(N)^G$ is surjective, as desired.

Moreover, if $\QQQ(\Gamma)=0$, then by \eqref{eq=mainexact} the map $i^{\ast}$ is also injective; hence $i^{\ast}$ is an isomorphism.
\end{proof}

Now, we proceed to the proof of  Theorem~\ref{thm=main}; it suffices to show that the group $\Gam$ of the form (a) or (b) fulfills conditions (1) and (2) of Proposition~\ref{prop=vanishingV}. For (1), the following is known.

\begin{theorem}\label{thm=Spc^2}
Let $\Gamma$ be a group of the form \textup{(a)} or \textup{(b)} as in Theorem~\textup{\ref{thm=main}}. Then, $\QQQ(\Gamma)=0$, and the second comparison map $c^2_{\Gam}\colon \HHH^2_b(\Gam;\RR)\to \HHH^2(\Gam;\RR)$ is an isomorphism.
\end{theorem}

We note that, under the condition $(\HHH^1_{/b}(\Gam)=)\QQQ(\Gamma)=0$, $c^2_{\Gam}$ is an isomorphism if and only if $c^2_{\Gam}$ is surjective; the surjectivity of $c^2_{\Gam}$ will be used when we apply Proposition~\ref{prop=vanishingV}.

\begin{proof}[Proof of Theorem~\textup{\ref{thm=Spc^2}}]
Note that $\Gam$ is a lattice in $\Sp(2g,\RR)$. The vanishing of $\QQQ(\Gamma)$ follows from \cite[Theorem~21]{BurgerMonod}. In what follows, we will see that $c^2_{\Gam}$ is an isomorphism. 

Consider the following commutative diagram:
\[
\xymatrix@C=5em{
    \HHH^2_{cb}(\Sp(2g,\RR);\RR) \ar[r]^-{c^2_{\Sp(2g,\RR)}}\ar[d]_{r_b} & \HHH^2_{c}(\Sp(2g,\RR);\RR) \ar[d]^{r}\\
     \HHH^2_{b}(\Gam;\RR) \ar[r]_{c^2_{\Gam}} & \HHH^2(\Gam;\RR). 
}
\]
Here, $c^2_{\Sp(2g,\RR)}$ is the second comparison map from countinuous bounded cohomology to continuous cohomology of $\Sp(2g,\RR)$; $r_b$ and $r$ are homomorphisms induced by the inclusion $\Gam \hookrightarrow \Sp(2g,\RR)$, respectively.

It is well known that the comparison map $c^2_{\Sp(2g,\RR)}$ of $\Sp(2g,\RR)$ gives isomorphisms $\HHH^2_{cb}(\Sp(2g,\RR);\RR)\cong \HHH^2_c(\Sp(2g,\RR);\RR)\cong \RR$; and that $r_b\colon \HHH^2_{cb}(\Sp(2g,\RR);\RR)\to \HHH^2_{b}(\Gam;\RR)$ is an isomorphism. We refer the reader to \cite[Section~6]{MS04} for details of these facts. Hence, it only remains to show that $r\colon \HHH^2_{c}(\Sp(2g,\RR);\RR)\to \HHH^2(\Gam;\RR)$ is an isomorphism. 

In case~(a), it is a classical fact that $r$ gives isomorphisms $\RR\cong \HHH^2_{c}(\Sp(2g,\RR);\RR)\cong \HHH^2(\Sp(2g,\ZZ);\RR)$ (see for instance \cite{BorelWallach}). In case~(b), Hain's result \cite[Theorem~3.2]{Hain}, which is related to the Borel stable range for $\mathrm{Sp}(2g,\mathbb{Z})$, shows that $r\colon \HHH^2_{c}(\Sp(2g,\RR);\RR)\to \HHH^2(\Gam;\RR)$ is an isomorphism. This completes our proof.
\end{proof}

Now, we go into the discussion on condition (2), which was the missing part of the proof of Theorem~\ref{thm=main} before Blatz's theorem (Theorem~\ref{thm=Blatz}). For the proof of Theorem~\ref{thm=SpH^3_b}, we also employ Monod's theorem \cite{Monod10}.

\begin{theorem}[corollary to the theorems of Blatz and Monod]\label{thm=SpH^3_b}
Let $g$ be an integer at least two. Let $\Gam$ be a lattice in $\Sp(2g,\RR)$. Then, we have $\HHH^3_b(\Gam;\RR)=0$. In particular, if $\Gam$ is a finite index subgroup of $\Sp(2g,\ZZ)$, then $\HHH^3_b(\Gam;\RR)=0$.
\end{theorem}

\begin{proof}
Observe that the (real) rank of $\Sp(2g,\RR)$ is $g$. Since $2g>3$,  Monod's theorem \cite[Corollary~1.4]{Monod10} applies to the lattice $\Gam$ in $\Sp(2g,\RR)$; hence, the inclusion map $\Gam\hookrightarrow \Sp(2g,\RR)$ induces an isomorphism $\HHH^3_{cb}(\Sp(2g,\RR);\RR)\cong \HHH^3_{b}(\Gam;\RR)$.  Blatz's theorem (Theorem~\ref{thm=Blatz}) implies that $\HHH^3_{cb}(\Sp(2g,\RR);\RR)=0$, and  this completes the proof. 
\end{proof}

\begin{proof}[Proof of Theorem~\textup{\ref{thm=main}}]
We can apply Proposition~\ref{prop=vanishingV} to $\Gam$ by Theorems~\ref{thm=Spc^2} and \ref{thm=SpH^3_b}.
\end{proof}

\begin{remark}
Our proof of Corollary~\ref{cor=CHH} is based on  Proposition~\ref{prop=vanishingV}, and hence is not constructive. It would be of independent interest to construct the extension explicitly, and to study the geometric meaning of it.
\end{remark}

\section{Proof of Theorem~\ref{thm=main_Johnson}}\label{section=proof2}
Recall that $\mathcal{I}_{g,1}$ denotes the Torelli group $\mathcal{I}(\Sigma_{g,1})$, and $\mathcal{K}_{g,1}$ denotes the Johnson kernel of $\Sigma_{g,1}$.
In this section, we use the symbol $\mathcal{M}_{g,1}$ to denote $\Mod(\Sigma_{g,1})$.

We use the following description of $\WW(G,N)$.
\begin{theorem}[{\cite[Theorem~1.10]{KKMMM}}]\label{thm=Wspace_isom}
  For an exact sequence $1 \to N \to G \xrightarrow{p} \Gam \to 1$ with boundedly $3$-acyclic $\Gam$, the space $\WW(G,N)$ is isomorphic to 
  \[
    \Im (p^* \colon \HHH^2(\Gam;\RR) \to \HHH^2(G;\RR)) \cap \Im (c_{G}^2 \colon \HHH_b^2(G;\RR) \to \HHH^2(G;\RR)).
  \]
\end{theorem}
Here a group $\Gam$ is said to be \emph{boundedly $n$-acyclic} if the bounded cohomology group $\HHH_b^k(\Gam;\RR)$ is trivial for $1 \leq k \leq n$.
It is known that abelian groups are boundedly $n$-acyclic for every $n$ (see \cite[Theorem 3.6]{Frigerio}).

Recall that Theorem~\ref{thm=main_Johnson} concerns $\WW(G,N)$ of the exact sequence \eqref{eq=exactJohnson}:
\[
  1\longrightarrow \mathcal{K}_{g,1} \longrightarrow \mathcal{I}_{g,1}\stackrel{\tau}{\longrightarrow} \Lambda^3H \longrightarrow 1.
\]
Since $\Lambda^3H$ is abelian and hence boundedly 3-acyclic, we have
\[
  \WW(\mathcal{I}_{g,1}, \mathcal{K}_{g,1}) \cong \Im (\tau^*) \cap \Im (c_{\mathcal{I}_{g,1}}^2)
\]
by Theorem~\ref{thm=Wspace_isom}.
Hence to prove Theorem \ref{thm=main_Johnson}, it suffices to construct a non-zero element of $\Im (\tau^*) \cap \Im (c_{\mathcal{I}_{g,1}}^2)$.

Let $\nu \colon \mathcal{M}_{g,1} \to \Homeo_+(S^1)$ be the Nielsen action and $\nu_I \colon \mathcal{I}_{g,1} \to \Homeo_+(S^1)$ the restriction of $\nu$ to $\mathcal{I}_{g,1}$.
Then the pullback of the Euler class $e \in \HHH^2(\Homeo_+(S^1);\RR)$ under $\nu_I$ is a nonzero element.
Indeed, the point-pushing subgroup $\pi_1(\Sigma_g) \subset \mathcal{I}_{g,1}$ acts on the Gromov boundary by a Fuchsian action, and the pullback of the Euler class under the Fuchsian action is nonzero (see \cite[Section 5.5]{FM12}).
Moreover, $\nu_I^* e$ is contained in $\Im (c_{\mathcal{I}_{g,1}}^2)$ since the Euler class $e$ is bounded.
Hence $\nu_I^* e$ is a nonzero element of $\Im (c_{\mathcal{I}_{g,1}}^2)$.

The Euler class $\nu_I^*e$ is also contained in $\Im (\tau^*)$.
This can be verified as follows.
Let $\rho_1 \colon \mathcal{M}_{g,1} \to \frac{1}{2} \Lambda^3H \rtimes \Sp(2g,\ZZ)$ be the extended Johnson homomorphism (\cite{MR1193604}), which gives rise to the commutative diagram
\[
\xymatrix{
  \mathcal{M}_{g,1} \ar[r]^-{\rho_1} & \frac{1}{2} \Lambda^3H \rtimes \Sp(2g,\ZZ) \\
  \mathcal{I}_{g,1} \ar@{_{(}-_>}[u]_-{i} \ar[r]^-{\tau} & \Lambda^3H \ar@{_{(}-_>}[u]_-{j}
}
\]
Morita proved in \cite{MR1659679} that the Euler class of $\mathcal{M}_{g,1}$ is contained in the image of the map $\rho_1^* \colon \HHH^*(\frac{1}{2} \Lambda^3H \rtimes \Sp(2g,\ZZ);\QQ) \to \HHH^*(\mathcal{M}_{g,1};\QQ)$  (see also \cite[Theorem 1]{MR1418577}).
By the commutativity of the above diagram, the Euler class $\nu_I^*e$ on $\mathcal{I}_{g,1}$ is contained in the image of $\tau^*$. 

We have shown that the Euler class $\nu_I^*e$ is nonzero and is contained in the intersection $\Im (\tau^*) \cap \Im (c_{\mathcal{I}_{g,1}}^2)$, which implies that $\WW(\mathcal{I}_{g,1}, \mathcal{K}_{g,1})$ is nonzero.

\section*{Acknowledgement} 
The first-named author, the second-named author, the fourth-named author and the fifth-named author are partially supported by JSPS KAKENHI Grant Number JP25K06994, JP24K16921, JP23K12975, and JP26K06799, respectively.
The third-named author is partially supported by JSPS KAKENHI Grant Number JP23K12971 and Research Fellowship Promoting International Collaboration, The Mathematical Society of Japan.

\section*{Use of AI}
We used ChatGPT for a literature search, which led us to \cite{Blatz}; the authors then consulted the original thesis and verified the statement and proof of Theorem 1.4 themselves.

\bibliography{reference}
\bibliographystyle{abbrv}

\end{document}